\documentclass[11pt]{article}
\usepackage[latin1]{inputenc}
\usepackage{amssymb, amsmath, amsthm, bm, epsfig}
\usepackage{graphicx, color}
\usepackage{pdflscape}

\newcommand {\tr} {\mbox{tr}}

\newcommand {\J}{\widetilde{\!\!\widetilde{J}}}

\title{Modified likelihood ratio tests in heteroskedastic multivariate regression models with 
measurement error}
\author
{
   Tatiane F. N. Melo\\
     {\it \footnotesize Instituto de Matemática e Estatística,
     Universidade Federal de Goiás}
     \vspace{-0.2cm}\\
     {\it \footnotesize Campus Samambaia, CP 131, Goiânia/GO, 74001-970, Brazil}
     \vspace{-0.2cm}\\
     {\footnotesize email: {\tt tmelo@mat.ufg.br}} \\ \\
   Silvia L. P. Ferrari\\
     {\it \footnotesize Departamento de Estatística, IME,
     Universidade de São Paulo}
     \vspace{-0.2cm}\\
     {\it \footnotesize Rua do Matão, 1010, São Paulo/SP, 05508-090, Brazil}
     \vspace{-0.2cm}\\
     {\footnotesize email: {\tt silviaferrari@usp.br}} \\ \\
  Alexandre G.~Patriota\\
     {\it \footnotesize Departamento de Estatística, IME,
     Universidade de São Paulo}
     \vspace{-0.2cm}\\
     {\it \footnotesize Rua do Matão, 1010, São Paulo/SP, 05508-090, Brazil}
     \vspace{-0.2cm}\\
     {\footnotesize email: {\tt patriota@ime.usp.br}} \\ \\     
}
\date{}

\begin{document}
\maketitle


\begin{abstract}
In this paper, we develop modified versions of the likelihood ratio test for multivariate heteroskedastic errors-in-variables regression models. The error terms are allowed to follow a multivariate distribution in the elliptical class of distributions, which has the normal distribution as a special case. We derive the Skovgaard adjusted likelihood ratio statistics, which follow a chi-squared distribution with a high degree of accuracy. We conduct a simulation study and show that the proposed tests display superior finite sample behavior as compared to the standard likelihood ratio test. We illustrate the usefulness of our results in applied settings using a data set from the WHO MONICA Project on cardiovascular disease. 
\end{abstract}

\vspace{0.1cm}

\section{Introduction}\label{sec1}

In actual problems, it is not uncommon to observe variables subject to measurement errors. This characteristic  leads to extra variations not predicted by the standard model. Such extra variations must be taken into account in statistical modeling to avoid misleading conclusions (e.g., inconsistent estimators, false confidence intervals, and underestimated standard errors). Fuller \cite{Fuller1987} presents many interesting problems with this peculiarity, one of his examples is presented in the following. Suppose that we wish to fit a regression line for corn yield  against nitrogen in the soil. It is known that the true quantity of nitrogen in the soil is not directly available and it has to be estimated by laboratory analysis from a selected sample. {In order to capture this extra variability of the collected data, measurement error equations are included in the modeling process. This procedure produces, in general, consistent estimators for the model parameters}. {Therefore, the} full specification of errors-in-variables regression models is essentially based in two sets of equations, one {related to the regression modeling} {(known as the regression equation) and other} for the measurement errors {(known as  measurement equation)}. 

{Measurement error models are widely applied in many fields, such as epidemiology, astronomy, biology, and so on. In the following we present some applications that motivated us the development of the present article.} Aoki et al. \cite{Aokietal2001} employ a statistical analysis to compare the effectiveness of two types of toothbrushes in removing plaque; the covariate (independent variable) is a plaque score before brushing and the response variable (dependent variable) is the plaque score after brushing. The amount of plaque before and after brushing is obtained in the same way. In this case, {as the amount of plaque is imprecisely evaluated, it is reasonable to consider that the covariate is subject to measurement errors}. Kulathinal et al. \cite{Kulathinaletal2002} and Patriota et al. \cite{Patriotaetal2009} study a univariate heteroskedastic errors-in-variables  model with equation error and normal distribution. An extension to the scale-mixture of normal distributions is presented by Cao et al. \cite{Caoetal2012}. The results in these three latter papers were applied to an epidemiological dataset of the WHO MONICA (World Health Organization Multinational MONItoring of trends and determinants in CArdiovascular disease) Project. Patriota et al. \cite{Patriotaetal2011} extend this univariate version to the multivariate one keeping the normality assumption, and also develop bias correction for the maximum likelihood estimators.  

In the aforementioned problems, the appropriate model should consider the measurement equation in the modeling, as  presented in Section \ref{sec2}. When the true covariates (the unobservable covariates) are modeled as random variables, the model is said to be structural. Otherwise, when these true covariates are fixed unknown parameters, the model is said to be functional. The structural model has identifiability problems while the functional version has unbounded likelihood. To overcome such problems, it is usual to regard some quantities as known, e.g., the measurement errors variances (\cite{ChanMak1979}, \cite{Wong1989}) or the model intercept \cite{Aokietal2001}. Here, we adopt the structural version; for more details regarding errors-in-variables models we refer the reader to Fuller \cite{Fuller1987}, Cheng and Van Ness \cite{ChengVanNess1999} and Buonaccorsi \cite{Buonaccorsi2010}.

Generally speaking, inference procedures under complicated models are based on asymptotic arguments. However, as shown by Melo and Ferrari \cite{MeloFerrari2010}, first-order asymptotics under errors-in-variables models may be inaccurate for small- and moderate-sized samples. With the purpose of improving the approximation of the distribution of the likelihood ratio statistic to its asymptotic reference chi-square distribution, Skovgaard \cite{Skovgaard2001} proposes a new test statistic when the interest lies in testing a parameter vector (see Section 3 for details). In this paper, we apply his proposal to a multivariate  regression model with measurement errors when the random quantities have a joint elliptical distribution. The model considered here is a generalization of the one studied by Patriota et al. \cite{Patriotaetal2011}. Therefore, our results apply to this latter paper and to all its subcases.

It is well-known that the normality assumption is not always tenable and alternative distributions (or methodologies) should be considered in such situations. For this reason, we consider the elliptical family of distributions, since this class includes many important distributions such as normal, Student-$t$, power exponential, contaminated normal, Pearson II, Pearson VII, and logistic, with heavier or lighter tails than the normal distribution.  Basically, we say that a $q\times 1$ random vector $\bm{Z}$ has a multivariate elliptical distribution with location parameter $\bm{\mu}$ ($q\times 1$) and a positive definite scale matrix $\bm{\Omega}$ ($q\times q$) if its density function exists, and is given by Fang et al. \cite{Fangetal1990}
\begin{equation}\label{dens}
f_{\bm{Z}}(\bm{z}) = |{\bm{\Omega}}|^{-1/2}
p_0\bigl[(\bm{z} - \bm{\mu})^{\top}{\bm{\Omega}}^{-1}(\bm{z} - \bm{\mu})\bigr],
\quad \bm{z}\in\mathbb{R}^q,
\end{equation} where $p_0:\mathbb{R}\to[0,\infty)$ is such that $\int_{0}^{\infty}u^{\frac{q}{2}-1}p_0(u) du<\infty$.  We use the notation $\bm{Z}\sim El_q(\bm{\mu}, {\bm{\Omega}}, p_0)$. It is possible to show that the characteristic function of $\bm Z$ is given by $\psi(t) = E(\exp(i\bm{t}^\top \bm{Z})) = \exp(i\bm{t}^\top \bm{\mu})\phi(\bm{t}^\top \bm{\Omega} \bm{t})$, where $\bm{t} \in \mathbb{R}^q$ and $\phi: [0,\infty) \to \mathbb{R}$. If $\phi$ is twice differentiable at zero, then $E(\bm{Z}) = \bm{\mu}$ and $\mbox{Var}(\bm{Z}) = c \bm{\Omega}$, with $c = -2\phi'(0)$. For more details, we refer the reader to Fang et al. \cite{Fangetal1990} and Fang and Anderson \cite{FangAnderson1990}.

The paper unfolds as follows. Section \ref{sec2} introduces the heteroskedastic multivariate regression model with measurement error. Section \ref{sec3} contains our main
results, namely the ancillary statistic and an explicit formula for the modified likelihood ratio test. Section \ref{sec4} presents a simulation study on the finite sample behavior of the standard likelihood ratio test and its modified counterpart. Our simulation results show that the likelihood ratio test tends to be liberal and its modified version is much less size-distorted. An application that uses real data is presented and discussed in Section \ref{sec5}. Finally, Section \ref{sec6} concludes the paper. Technical details are collected in three appendices.  

\section{Heteroskedastic multivariate regression model with measurement error}
\label{sec2}

We consider the following heteroskedastic multivariate regression model with measurement error
\begin{equation}\label{E.4.1}
\begin{split}
\bm{y}_i &= \bm{\alpha} + \bm{\beta} {\bm x}_i + \bm{q}_i, 
\\
\bm{Y}_i &= \bm{y}_i + \bm{e}_i,
\\
{\bm X}_i &= {\bm x}_i + {\bm u}_i,
\end{split}
\end{equation}
for $i = 1, 2, \ldots, n$, where $\bm{y}_i$ is an ($m \times 1$) latent  response vector, $\bm{x}_i$ is a ($p \times 1$) latent vector of covariates, $\bm{\alpha}$ is a vector of intercepts, and $\bm{\beta}$ is an $(m \times p)$ matrix of parameters. The error terms are $\bm{q}_i, \bm{e}_i$ and $\bm{u}_i$. The vectors $\bm{y}_i$, and $\bm{x}_i$ are not directly observed; instead, the vectors $\bm{Y}_i$ and $\bm{X}_i$, respectively, are observed. The observed and unobserved vectors are connected by the additive relation specified by the two last measurement equations in (\ref{E.4.1}). Model (\ref{E.4.1}) can be shortly written as 
\begin{equation}\label{E.4.3}
{\bm Z}_i = {\bm\delta} + \Delta {\bm b}_i,
\end{equation}
with $i = 1, 2, \ldots, n$, where
\begin{equation*}\label{E.4.4}
{\bm Z}_i = 
\begin{pmatrix}
\bm{Y}_i \\ 
{\bm X}_i 
\end{pmatrix}, \quad  
{\bm\delta} = \begin{pmatrix}\bm{\alpha} \\ \bm{0} \end{pmatrix}, \quad 
\Delta = \begin{pmatrix} 
{\bm \beta}  & I_m         & I_m           & {\bm 0} \\ 
I_p          &{\bm 0}      & {\bm 0}       &  I_p 
\end{pmatrix}, \quad {\rm and} \quad 
{\bm b}_i = 
\begin{pmatrix} {\bm x}_i \\ \bm{q}_i \\ \bm{e}_i \\{\bm u}_i \end{pmatrix},
\end{equation*}
where $I_p$ is the identity matrix of dimension $p\times p$. We assume that ${\bm b}_1, {\bm b}_2, \ldots, {\bm b}_n$ are independent and identically distributed according to an elliptical distribution with the function $\phi$ being twice differentiable at zero (see the previous section, below Eq.~(\ref{dens})). This model specification generalizes the one proposed by Patriota et al. \cite{Patriotaetal2011}. Here, we define the covariance matrix of $\bm{b}_i$ as
\[ \mbox{Var}(\bm{b}_i) = \Pi_i = 
 \begin{pmatrix}
\Sigma_x   & {\bm 0}    & {\bm 0}         & {\bm 0} \\ 
{\bm 0}    & \Sigma_{q} & {\bm 0}         & {\bm 0} \\ 
{\bm 0}    & {\bm 0}    & \Sigma_{e_i}    & \Sigma_{(ue)_i}^\top \\ 
{\bm 0}    & {\bm 0}    & \Sigma_{(ue)_i} & \Sigma_{u_i}
\end{pmatrix},
\] 
where $\Sigma_{e_i}$, $\Sigma_{(ue)_i}$, and $\Sigma_{u_i}$ are known $m\times m$, {$p \times m$} and $p \times p$ matrices, respectively. If these quantities do not depend on the index $i$ then we have the homoskedastic multivariate regression model with measurement error. In short 
${\bm b}_i \sim El_{(2m+2p)}(\bm{\eta}, \Omega_{\bm{b}_i}; p_0)$, with
\begin{eqnarray*}\label{E.4.5}
{\bm \eta} = \begin{pmatrix} {\bm \mu}_x \\ \bm{0} \\ {\bm 0} \\ {\bm 0}\end{pmatrix}, \quad
{\Omega}_{\bm{b}_i} = c^{-1}\Pi_i,
\end{eqnarray*}
and $c = -2\phi'(0)$. It can be shown that the vector ${\bm Z}_i$ has a $(p+m)$-variate elliptical distribution with location vector ${\bm \mu} = {{\bm \mu}({\bm \theta})} \in \mathbb{R}^{p+m}$ and dispersion matrix $\Omega_i = \Omega_i({\bm \theta})$. The parameter vector is ${\bm \theta} = (\mbox{vec}(\bm {\beta})^\top, \bm{\alpha}^\top, {\bm \mu}_x^\top, \mbox{vech}(\Sigma_q)^{\top},$ $ \mbox{vech}(\Sigma_x)^\top)^\top \in \Theta \subset \mathbb{R}^s$, where $s =  mp + m + p +p(p+1)/2 + m(m+1)/2$ {is the dimension of the parameter space}, vec is the vector operator stacking the columns of a matrix underneath each other, and vech is the operator that stacks the lower triangle of a matrix into a column vector. 
{By using properties of elliptical distributions, we} can write  ${\bm Z}_i \sim El_{p+m}({\bm \mu}, \Omega_i; p_0)$, where ${\bm \mu} = {\bm \delta} + \Delta{\bm \eta}$ and $\Omega_i = \Delta {\Omega}_{\bm{b}_i} \Delta^\top$ {are the location vector and scale matrix, respectively, which are functions of the parameter vector $\bm{\theta}$, namely:}
\begin{eqnarray*}\label{E.4.6}
{{\bm \mu}} {\equiv \bm{\mu}(\bm{\theta})=} 
\begin{pmatrix}
\bm{\alpha} + {\bm \beta} {\bm \mu}_x \\ {\bm \mu}_x 
\end{pmatrix} \quad {\rm and}\quad
\Omega_i {\equiv \Omega_i(\bm{\theta})=}  
c^{-1}\begin{pmatrix}
{\bm \beta} \Sigma_x {\bm \beta}^\top + \Sigma_q + \Sigma_{e_i} & {\bm \beta}\Sigma_x + \Sigma_{(ue)_i}^\top\\ 
\Sigma_x {\bm \beta}^\top + \Sigma_{(ue)_i}                     & \Sigma_x + \Sigma_{u_i}
\end{pmatrix}.
\end{eqnarray*}

{Notice that, for the normal case,  model (\ref{E.4.3}) is identifiable. Hence, if the density generator $p_0$ does not depend on additional parameters to be estimated, model (\ref{E.4.3}) is also identifiable for the class of elliptical distributions. It is noteworthy that for many elliptical distributions $p_0$ contains unknown quantities, e.g., the degrees of freedom in the Student-$t$ distribution and the shape parameter in the power exponential distribution. We may want to estimate these quantities via maximum-likelihood estimation, however, Lucas  \cite{Lucas1997} studies some robust features of M-estimators by using influence functions for the Student-$t$ distribution and shows that the protection against ``large'' observations is only valid when the degrees of freedom parameter is known. Therefore, for the purpose of avoiding possible lack of protection against outliers, we do not estimate those unknown parameters. Otherwise, unboundedness problems may arise for the influence functions and the elliptical distribution will lose its main goal. In practice, we can use model selection procedures, such as the Akaike information criterion (AIC), to choose the more appropriate values of such unknown parameters.} 

The likelihood function for model (\ref{E.4.3}) is 
\begin{equation*}\label{E.4.7}
L(\theta) = p_{\bm Z}({\bm z}, {\bm \theta}) =  \prod_{i=1}^{n} |\Omega_i|^{-1/2}
p_0({\bm d}_i^\top \Omega_i^{-1} {\bm d}_i),
\end{equation*}
where ${\bm d}_i = {\bm d}_i({\bm \theta}) = {\bm z}_i - {\bm \mu}$. 
The logarithm of the likelihood function is
\begin{eqnarray}\label{E.4.8}
\ell({\bm \theta}) = \sum_{i=1}^{n} \left[-\frac{1}{2} \log |\Omega_i| + \log p_0
({\bm d}_i^\top \Omega_i^{-1} {\bm d}_i)\right].
\end{eqnarray}

Lemonte and Patriota \cite{LemontePatriota2011} discuss the maximum likelihood estimator of $\bm \theta$, and derive the generalized leverage as well as the normal curvatures of local influence under some perturbation schemes. {The reader is referred to their work for further details on the iterative procedure to attain the maximum likelihood estimates}.

\section{Modified likelihood ratio test}
\label{sec3}

Let ${\bm \theta} = ({\bm \psi}^\top, {\bm \omega}^\top)^\top$, where ${\bm \psi}$ is the vector of parameters of interest and ${\bm \omega}$ is the vector of nuisance parameters. Here, the null and alternative hypotheses are, respectively: ${\cal H}_{0}: {{\bm \psi}} = {{\bm \psi}}^{(0)}$ and ${\cal H}_{1}: {{\bm \psi}} \neq {{\bm \psi}}^{(0)}$, where  ${{\bm \psi}}^{(0)}$ is a known $q$-vector. In this paper, we focus on testing hypotheses on some elements of the matrix $\bm{\beta}$. Define $\mbox{vec}(\bm \beta) = (\bm{\beta}_1^\top, \bm{\beta}_2^\top)^\top$; without loss of generality,  we consider ${\bm \psi} = \bm{\beta}_1$ and  ${\bm \omega} = (\bm{\beta}_2^\top, \bm{\alpha}^\top, {\bm \mu}_x^\top, \mbox{vech}(\Sigma_q)^\top, \mbox{vech}(\Sigma_x)^\top)^\top$.  The maximum likelihood estimator over $\Theta$ is denoted by $\bm{\widehat{\theta}} = (\bm{\widehat{\psi}}, \bm{\widehat{\omega}}^\top)^\top$  and the maximum likelihood estimator under the null hypothesis, by $\bm{\widetilde\theta} = (\bm{\widetilde\psi}, \bm{\widetilde\omega} ^\top)^\top$, where $\bm{\widetilde\psi} = 
{\bm{\psi}}^{(0)}$. We use ``\ $\widehat{}$\ " and ``\ $\widetilde{}$\ " for matrices and vectors to indicate that they are computed at $\bm{\widehat{\theta}}$ and $\bm{\widetilde\theta}$, respectively.

The likelihood ratio statistic for testing ${\cal H}_0$ is
\begin{eqnarray*}\label{E.4.11}
LR = 2\:\left\{\ell(\bm{\widehat{\psi}}) - \ell(\bm{\widetilde{\psi}})\right\}.
\end{eqnarray*}
Under regularity conditions on $p_0$ see, for instance, Cox and Hinkley \cite{CoxHinkley1974}, this statistic converges in distribution, under ${\cal H}_{0}$, to ${\cal X}^2_q$, where $q$ is the dimension of ${\bm \psi}$. 

When the parameter of interest, $\bm{\psi}$, is one-dimensional, we can use a signed likelihood ratio statistic to test ${\cal H}_{0}$ against ${\cal H}_{1}$. Barndorff-Nielsen \cite{Barndorff1991} proposes a modified version of this statistic that depends on an ancillary statistic and some derivatives with respect to the sample space. This modified version intends to better approximate the signed likelihood ratio statistic distribution for the standard normal distribution.
The extension for the vectorial case was developed by Skovgaard \cite{Skovgaard2001}, for which two modified statistics, that are asymptotically equivalent, were proposed. We denote these two statistics here by $LR_a^*$ and $LR_a^{**}$. These modified statistics retain the essential character of $LR$, but can be difficult to obtain. The difficulty arises from the need of an appropriate ancillary statistic and derivatives of the log-likelihood function with respect to the data.  By ``ancillary statistic'' we mean a statistic, say ${\bm a}$, whose distribution does not depend on the unknown parameter $\bm{\theta}$, and that along with the maximum likelihood estimator $\bm{\widehat{\theta}}$, is a minimal sufficient statistic for the model. If $(\bm{\widehat{\theta}},{\bm a})$ is sufficient, but not minimal sufficient, Barndorff-Nielsen's results still hold; see, Severini \cite{Severini2000} (\S \ 6.5). In fact, minimal sufficiency is only required for the ancillary ${\bm a}$ to be relevant to the statistical analysis. Sufficiency implies that the log-likelihood function depends on the data only through $(\bm{\widehat{\theta}}, {\bm a})$, and we then write $\ell(\bm{\theta};\bm{\widehat \theta},{\bm a})$. The required derivatives of  $\ell(\bm{\theta};\bm{\widehat \theta},{\bm a})$ with respect to the data are
\begin{eqnarray*}\label{E.17}
\ell' = \frac{\partial \ell(\bm{\theta};\bm{\widehat \theta},a)}
{\partial \bm{\widehat \theta}}, \ \  
U'= \frac{\partial^2 \ell(\bm{\theta};\bm{\widehat \theta},
{\bm a})}{\partial \bm{\widehat \theta} \partial{\bm{\theta}}^\top} \ \ {\rm and} \ \ 
{\:\:{\J}}= \frac{\partial^2 \ell(\bm{\theta};\bm{\widehat \theta},{\bm a})}{\partial \bm{\widehat \theta} \partial{\bm{\theta}^\top}}\Bigg|_{\bm{\widehat \theta} = \bm{\widetilde \theta},\:\: \bm{\theta} = \bm{\widetilde \theta}}.
\end{eqnarray*}

The modified likelihood ratio statistics are
\begin{eqnarray*}\label{E.4.15}
LR_a^* = LR \left(1 - \frac{1}{LR} \log\rho \right)^2
\end{eqnarray*}
and
\begin{eqnarray*}\label{E.4.16}
LR_a^{**} = LR - 2\log\rho,
\end{eqnarray*}
with
\begin{eqnarray}\label{E.4.17}
\rho = |\widehat{J}\:|^{1/2} |{\widetilde U}'|^{-1} |
{\widetilde J}_{\bm{\omega\omega}}|^{1/2} |
\:{\J}_{\bm{\omega\omega}}|^{-1/2} |\:{\J}\:|^{1/2} 
\frac{\{{\widetilde U}^{\top} {\J}^{\: -1} {\widetilde U}\}^{p/2}}
{LR^{q/2 - 1} ({\widehat \ell}'- {\widetilde \ell}')^{\top} 
({\widetilde U}')^{-1} {\widetilde U}},
\end{eqnarray}
The quantities ${\widehat \ell}'- {\widetilde \ell}'$, ${\widetilde U}'$ and $\:\:{\J} = J(\bm{\widetilde{\theta}}; \bm{\widetilde{\theta}}, {\bm a})$ are derivatives with respect to the sample space and are obtained from an ancillary statistic ${\bm a}$ and the maximum likelihood estimator $\bm{\widehat{\theta}}$. Matrix ${\:\:{\J}}_{\bm{\omega\omega}}$ is the lower-right sub-matrix of ${\:\:{\J}}$ related to the nuisance parameters $\bm{\omega}$. Although the statistic $LR_a^*$ is non-negative,  and when $q=1$, it reduces to the square of Barndorff-Nielsen's statistic, the second version, $LR_a^{**}$, seems to be numerically more stable  and is naturally attained from theoretical developments see Skovgaard \cite{Skovgaard2001}. These statistics can be approximated to ${\cal X}^2_q$ with high accuracy under the null hypothesis.

In order to find an ancillary statistic, we first note that model (\ref{E.4.3}) is a transformation model. Hence any maximal invariant statistic 
is an ancillary statistic; see Barndorff-Nielsen et al \cite{Barndorffetal1989} (Chapter 8) and Barndorff-Nielsen \cite{Barndorff1986}. Moreover, a statistic is maximal invariant if it is a function of any invariant statistic. Pace and Salvan \cite{PaceSalvan1997} (Theorem 7.2) showed that in transformation models all invariant statistics are distribution constant, i.e., their distributions do not depend on $\bm{\theta}$. {Let $P_i\equiv P_i(\theta)$ be a lower triangular matrix such that $P_i P_i^\top = \Omega_i$ is the Cholesky decomposition for all $i=1, \ldots, n$.} By using the same idea as Melo and Ferrari \cite{MeloFerrari2010}, the statistic  ${\bm a} = ({\bm a}_1^\top, {\bm a}_2^\top$, $\ldots, {\bm a}_n^\top)^\top$, with
\begin{equation}\label{E.20}
{\bm a}_i = {\widehat P}_i^{-1}\left({\bm z}_i - {\bm{\widehat \mu}}\right),
\end{equation}
is maximal invariant and consequently is ancillary, as previously shown, { where $\widehat{\bm{\mu}} = \bm{\mu}(\widehat{\bm{\theta}})$ and $\widehat{P}_i = P_i(\widehat{\bm{\theta}})$} are the maximum likelihood estimates of ${\bm{\mu}}$ and ${P}_i$, respectively. Taking this ancillary statistic, we can find the derivatives with respect to the sample space required for Skovgaard's adjustment.

In order to compute the score vector, observed Fisher information, and the derivatives with respect to the sample space, we define the following matrices and vectors. Let $T$ be a matrix, with its $(j,k)$-th element given as
\begin{equation*}\label{E.4.18}
t_{jk} = \sum_{i=1}^{n} \tr\left({\Omega}_i^{-1} \Omega_{\theta_j\theta_k} - \Omega_i^{-1} {\Omega}_{\theta_k} \Omega_i^{-1} {\Omega}_{\theta_j}\right),
\end{equation*}  
where ${\Omega}_{\theta_j} = {\partial\Omega_i}/{\partial \theta_j}$ and $\Omega_{\theta_j\theta_k} = {\partial{\Omega}_{\theta_j}}/{\partial \theta_k}$, for $j, k = 1, 2, \ldots, s$. Observe that ${\Omega}_{\theta_j}$ does not depend on the index $i$, since the matrix $\Omega_i$ merely depends on $i$ through the known matrices $\Sigma_{e_i}$, $\Sigma_{(ue)_i}$, and $\Sigma_{u_i}$. Therefore, their derivatives with respect to the model parameters are zero. Consider the following block-diagonal matrices  $R = {\rm diag}({\bm r}, {\bm r}, \ldots, {\bm r})$ and $V = {\rm diag}({\bm v}, {\bm v}, \ldots, {\bm v})$, of dimension  $sn \times s$. The $i$-th elements of ${\bm r}$ and ${\bm v}$ are, respectively, ${r}_i = W_{p_0}({\bm d}_i^\top \Omega_i^{-1} {\bm d}_i)$ and ${v}_i = W_{p_0}'({\bm d}_i^\top \Omega_i^{-1} {\bm d}_i)$. 
Let ${\bm h} = \left({{\bm h}^{(1)^\top}}, \ldots, {{\bm h}^{(s)^\top}}\right)^\top$ and ${\bm w} = \left({{\bm w}^{(1)^\top}}, \ldots, {{\bm w}^{(s)^\top}}\right)^\top$ be $sn$-vectors, for which the $i$-th elements are given, respectively, by ${h}_i^{(j)} = -{\bm d}_i^\top \Omega_i^{-1} {\Omega}_{\theta_j} \Omega_i^{-1} {\bm d}_i - 2 \bm{\mu}_{\theta_j}^\top  {\Omega}_i^{-1} {\bm d}_i$, and ${w}_i^{(j)} = \left({\widehat P}_{\theta_j} {\bm a}_i + \bm{\widehat \mu}_{\theta_j}\right)^\top \Omega_i^{-1} {\bm d}_i$, where $P_{\theta_j} = {\partial P_i}/{\partial \theta_j}$ and $\bm{\mu}_{\theta_j} = {\partial\bm{\mu}}/{\partial \theta_j}$. In addition, let $B$, $C$, $F$, $G$, $M$, and $Q$ be ($sn \times s$) block-matrices, whose $(j,k)$-th blocks are, respectively, given by the vectors ${\bm b}^{(jk)}$, ${\bm c}^{(jk)}$, ${\bm f}^{(jk)}$, ${\bm g}^{(jk)}$, ${\bm m}^{(jk)}$, and ${\bm q}^{(jk)}$.  The $i$-th components of these vectors are, respectively, given by ${b}_i^{(jk)} = {h}_i^{(k)} {h}_i^{(j)}$, $c_i^{(jk)} = {w}_i^{(k)} {h}_i^{(j)}$,

\begin{equation*}\label{E.4.21}
\begin{split}
{f}_i^{(jk)} &= -\left({\widetilde P}_{\theta_k} {\bm a}_i + \bm{\widetilde \mu}_{\theta_k}
\right)^\top {\widetilde\Omega}_i^{-1} {{\widetilde\Omega}}_{\theta_j} {\widetilde\Omega}_i^{-1} {\widetilde P}_i {\bm a}_i - \bm{\widetilde\mu}_{
\theta_j}^\top {\widetilde\Omega}_i^{-1}\left({\widetilde P}_{\theta_k} {\bm a}_i + \bm{
\widetilde \mu}_{\theta_k}\right),
\\ \\
{g}_i^{(jk)} &= ({\widetilde P}_{\theta_k} {\bm a}_i + \bm{\widetilde \mu}_{\theta_k})^\top 
{\widetilde \Omega}_i^{-1} {\widetilde P}_i {\bm a}_i \left(-{\bm a}_i^\top {\widetilde P}_i^\top  
{\widetilde\Omega}_i^{-1} {{\widetilde\Omega}}_{\theta_j} {\widetilde\Omega}_i^{-1}{\widetilde P}_i {\bm a}_i - 2 \bm{\widetilde\mu}_{\theta_j}^\top 
{\widetilde\Omega}_i^{-1}{\widetilde P}_i {\bm a}_i\right),
\\ \\
{m}_i^{(jk)} &= {\bm d}_i^\top \left(2 {\Omega}_i^{-1} {\Omega}_{\theta_k} {\Omega}_i^{-1} {\Omega}_{\theta_j} {\Omega}_i^{-1} - {\Omega}_i^{-1}       \Omega_{\theta_j\theta_k}{\Omega}_i^{-1}\right) {\bm d}_i 
+ 2\bm{\mu}_{\theta_k}^\top {\Omega}_i^{-1} {\Omega}_{\theta_j} {\Omega}_i^{-1} {\bm d}_i 
+ 2\bm{\mu}_{\theta_j}^\top {\Omega}_i^{-1} {\Omega}_{\theta_k} {\Omega}_i^{-1} {\bm d}_i 
\\
& \ \ \ \ \ - 2 \bm{\mu}_{\theta_j\theta_k}^\top  {\Omega}_i^{-1} {\bm d}_i 
+ 2 \bm{\mu}_{\theta_j}^\top {\Omega}_i^{-1} \bm{\mu}_{\theta_k},
%
\\ \\
{q}_i^{(jk)} &= -\left({\widehat P}_{\theta_k} {\bm a}_i + \bm{\widehat \mu}_{\theta_k}
\right)^\top \Omega_i^{-1} {\Omega}_{\theta_j} \Omega_i^{-1} {\bm d}_i - \bm{\mu}_{\theta_j}^\top\Omega_i^{-1}\left({\widehat P}_{\theta_k} 
{\bm a}_i + \bm{\widehat \mu}_{\theta_k}\right),
\end{split}
\end{equation*}
where  $\bm{\mu}_{\theta_j\theta_k} = {\partial\bm{\mu}_{\theta_j}}/{\partial \theta_k}$. The score vector, 
observed information matrix, and derivatives with respect to the sample space are given in matrix notation by
\begin{equation*}\label{E.4.23}
U = - \frac{1}{2}\: {\bm n}^* + R^\top {\bm h}, \:\:\: J = \frac{1}{2}\:  T - R^\top M - 
V^\top Q, 
\end{equation*}
\begin{equation*}\label{E.4.24}
{\ell'} = 2 R^\top {\bm w}, \:\:\: {U'} = 2\left(R^\top B + V^\top C\right),\:\:\:{\rm and}\:\:\: 
{\:\:{\J}} = 2\left({\widehat R}^\top F + {\widehat V}^\top G\right),
\end{equation*}
where the $j$-th component of the vector ${\bm n}^*$ is $\sum_{i=1}^n \tr({\Omega}_i^{-1} {\Omega}_{\theta_j})$.

When we insert $\widehat J$, ${\widetilde J}_{\bm{\omega\omega}}$, ${\:{\J}}_{\bm{\omega\omega}}$, ${\:{\J}}$, $\widetilde U$, ${\widetilde U}'$, ${\widehat \ell}'- {\widetilde \ell}'$, and the original likelihood ratio statistic $LR$ into (\ref{E.4.17}), we obtain the required quantity $\rho$ for Skovgaard's adjustment. Now, one is able to compute the modified versions $LR_a^{*}$ and $LR_a^{**}$.

Computer packages that perform simple operations on matrices and vectors can be used to calculate $\rho$. Note that $\rho$ depends on the model through $\bm{\mu}$, $P_i$, $\Omega_i$ and $\Omega_i^{-1}$. The dependence on the specific distribution of ${\bm z}$ in the class of elliptical distributions occurs through $W_{p_0}$. Appendix A gives the required derivatives of ${\bm{\mu}}$ and ${\Omega}_i$ with respect to the vector of parameters. The derivative of $P_i$ with respect to the vector of parameters was obtained by using the algorithm proposed by Smith \cite{Smith1995}. 	           

\section{Simulation study}
\label{sec4}

In this section, we present simulation results to illustrate the effectiveness of the Skovgaard's adjustments for the original likelihood ratio statistic. The performance of the statistics $LR$, $LR_a^{*}$, and $LR_a^{**}$ were evaluated considering the repeated sampling principle. We study the frequencies of simulated samples that commit the type I error, and compare them with actual pre-assigned nominal levels.

The simulation study was based on model $(\ref{E.4.3})$ when $m=1$ and the random $(p+1)$-vector ${\bm Z}_i$ follows a multivariate normal, Student-$t$ with $\nu$ degrees of freedom or power exponential distribution with shape parameter $\lambda$ (Gómez et al. \cite{Gomezetal1998}). 
This study was conducted using the matrix programming language {\tt Ox} Doornik \cite{Doornik2006} considering $10,$$000$ (ten thousand) Monte Carlo samples. We consider the following sample sizes: $n = 20, 30, 40$ and $50$. The nominal levels are  $\gamma=1\%$, $5\%$, and $10\%$. In addition, we consider two to four covariates, i.e.,  $p = 2, 3$, and $4$. The null and alternative hypotheses are, respectively, ${\mathcal H}_{0}: {\bm \psi} = {\bm 0}$ and ${\mathcal H}_{1}: {\bm \psi} \neq {\bm 0}$, where ${\bm \psi} = (\beta_1, \beta_2, \ldots, \beta_q)^\top$, for $q = 2, 3$ or $4$.  The true parameters are $\alpha = 0.2$, ${\bm \mu}_x = (-2, -2, \ldots, -2)^\top$, $\Sigma_q = 10$ and $\Sigma_x = 4 I_p$, where  $I_p$ is the $(p \times p)$ identity matrix. As before $\Sigma_{e_i}$, $\Sigma_{(ue)_i}$ and $\Sigma_{u_i}$, for $i = 1, 2, \ldots, n$, are assumed to be known. We
set $\Sigma_{(ue)_i}$ as a $p \times 1$ vector of zeros, the values for $\Sigma_{e_i}$, for $i=1,\ldots, n$ were generated as $n$ independent observations from $\sqrt{\Sigma_{e_i}} \sim U(0,1)$,  and ${\Sigma_{u_i}}$, for $i=1,\ldots n$, was assumed to be a diagonal matrix; the square root of its elements were independently generated from a $U(0,1)$ distribution. {After generating these scale matrices, they were kept fixed for all Monte Carlo Simulations.} For the Student-$t$ distribution, we fixed the degrees of freedom at $\nu = 5$ and, for power exponential model, the shape parameter was fixed at 
$\lambda = 0.6$. Tables 1, 2 and 3 present the rejection rates for the tests based on $LR$, $LR_a^*$, and $LR_a^{**}$, for $q = 2$ and different values for $p$ and $n$.

In our simulations we observe that tests based on the original likelihood ratio statistic $LR$ are liberal when the sample size is small, since the respective rejection rates are larger than the actual nominal levels. This can be seen in Table 2 when $\bm{Z}_i$ follows a Student-$t$ distribution, $p = 3$, $q = 2$, and $n=20$.  In this case, the rejection rates for tests based on $LR$ equal $2.5\%$ $(\gamma=1\%)$, $9.0\%$ $(\gamma=5\%)$, and $15.7\%$ $(\gamma=10\%)$.  When ${\bm Z}_i$ follows a normal distribution, $p = 2 = q$, and $n=20$, we observe the following rejection rates: $2.1\%$ $(\gamma=1\%)$, $8.2\%$ $(\gamma=5\%)$, and $14.7\%$ $(\gamma=10\%)$; see Table 1.

On the other hand, tests based on the modified versions $LR_a^{*}$ and $LR_a^{**}$ present rejection rates closer to the actual nominal levels than the original version $LR$. Table 2 gives the rejection rates when $\bm{Z}_i$ follows a Student-$t$ distribution, $p = 4$, and $q = 2$. For $n=20$ and $\gamma=10\%$, the rejection rates are $16.1\%$ $(LR)$, $10.9\%$ $(LR_a^{*})$, and  $10.2\%$ $(LR_a^{**})$. For $n=30$ and $\gamma =10\%$, the rejection rates are  $13.9\%$ $(LR)$, $10.2\%$ $(LR_a^{*})$, and $10.0\%$ $(LR_a^{**})$. Another example can be seen in Table 1 for the normal distribution, $p = 3$, $q = 2$, $n = 30$, and $\gamma=1\%$. For this case, the rejection rates are $2.2\%$ $(LR)$, $1.2\%$ $(LR_a^{*})$, and $1.2\%$ $(LR_a^{**})$. For the power exponential distribution, $p = 4$, $q = 2$, $n = 30$, and $\gamma=5\%$ we have $7.8\%$ $(LR)$, $5.1\%$ $(LR_a^{*})$, and $4.9\%$ $(LR_a^{**})$; see Table 3.

\begin{figure}[!htp]  
\begin{center}
\includegraphics[width=16cm]{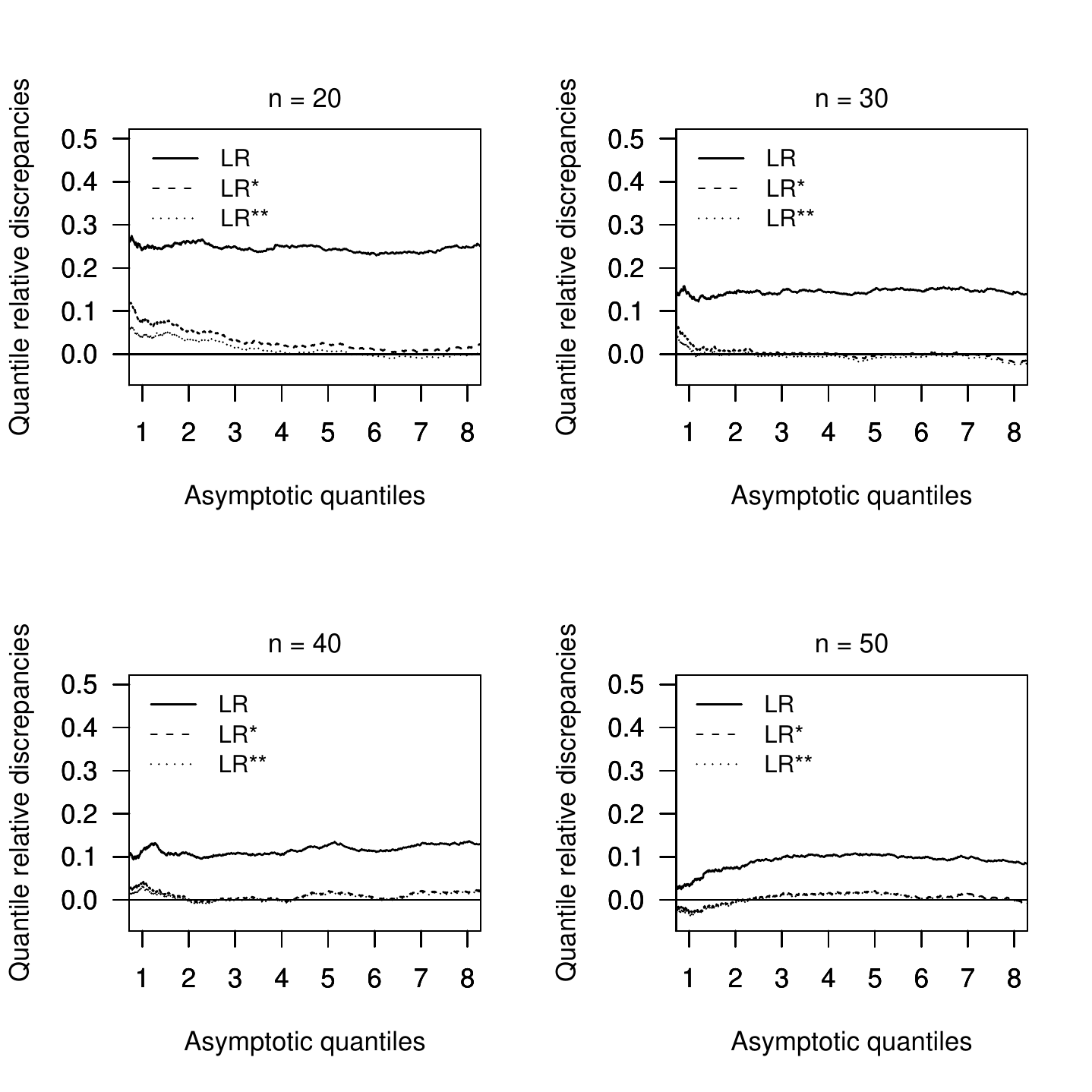}
\caption{Quantile relative discrepancies: $p = 4$, $q = 3$, and Student-$t$ distribution for ${\bm Z}_i$.}
\label{tStudp4q3}
\end{center}
\end{figure}

\begin{figure}[!htp]  
\begin{center}
\includegraphics[width=16cm]{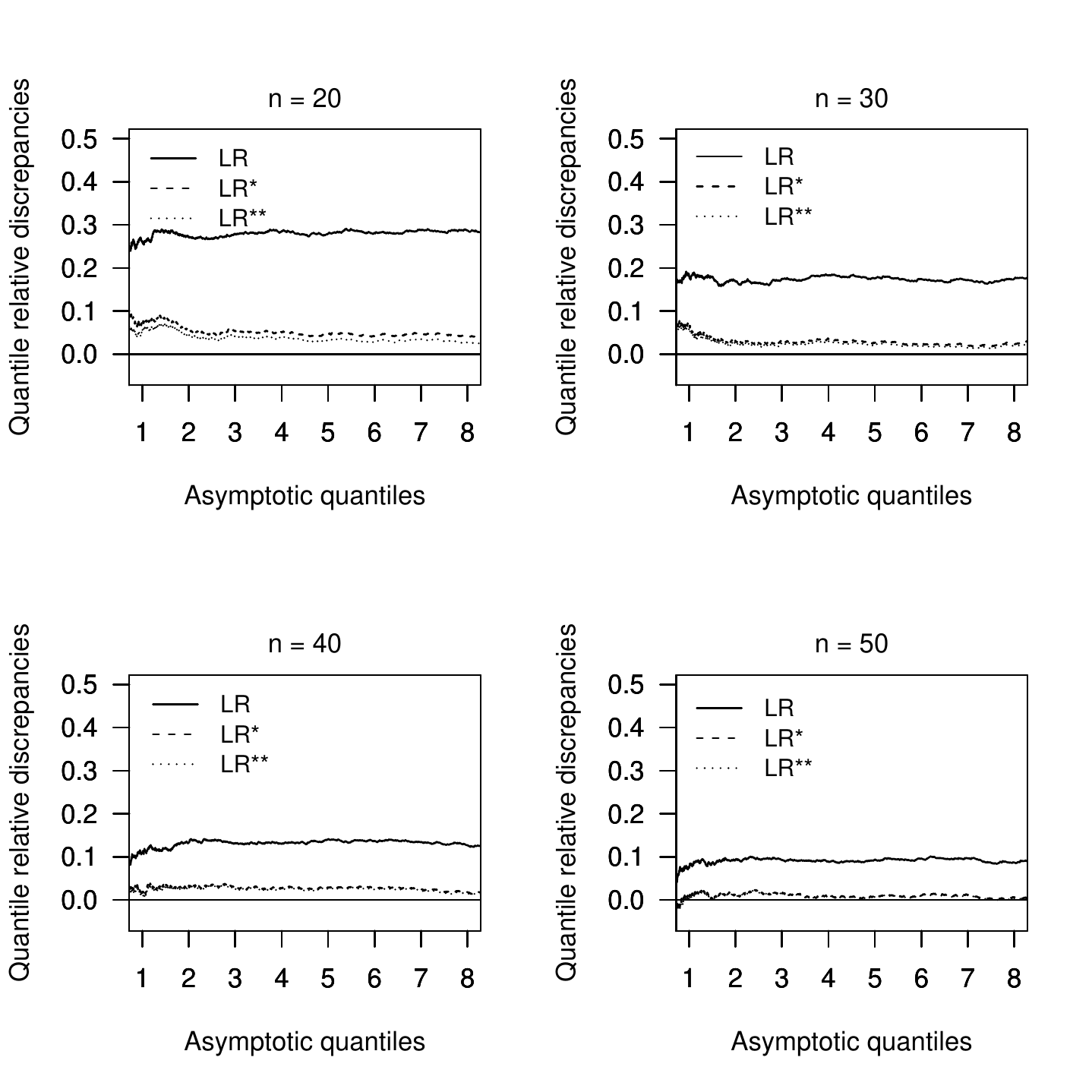}
\caption{Quantile relative discrepancies: $p = 4$, $q = 4$, and normal distribution for ${\bm Z}_i$.}
\label{Normp4q4}
\end{center}
\end{figure}

\begin{figure}[!htp]  
\begin{center}
\includegraphics[width=16cm]{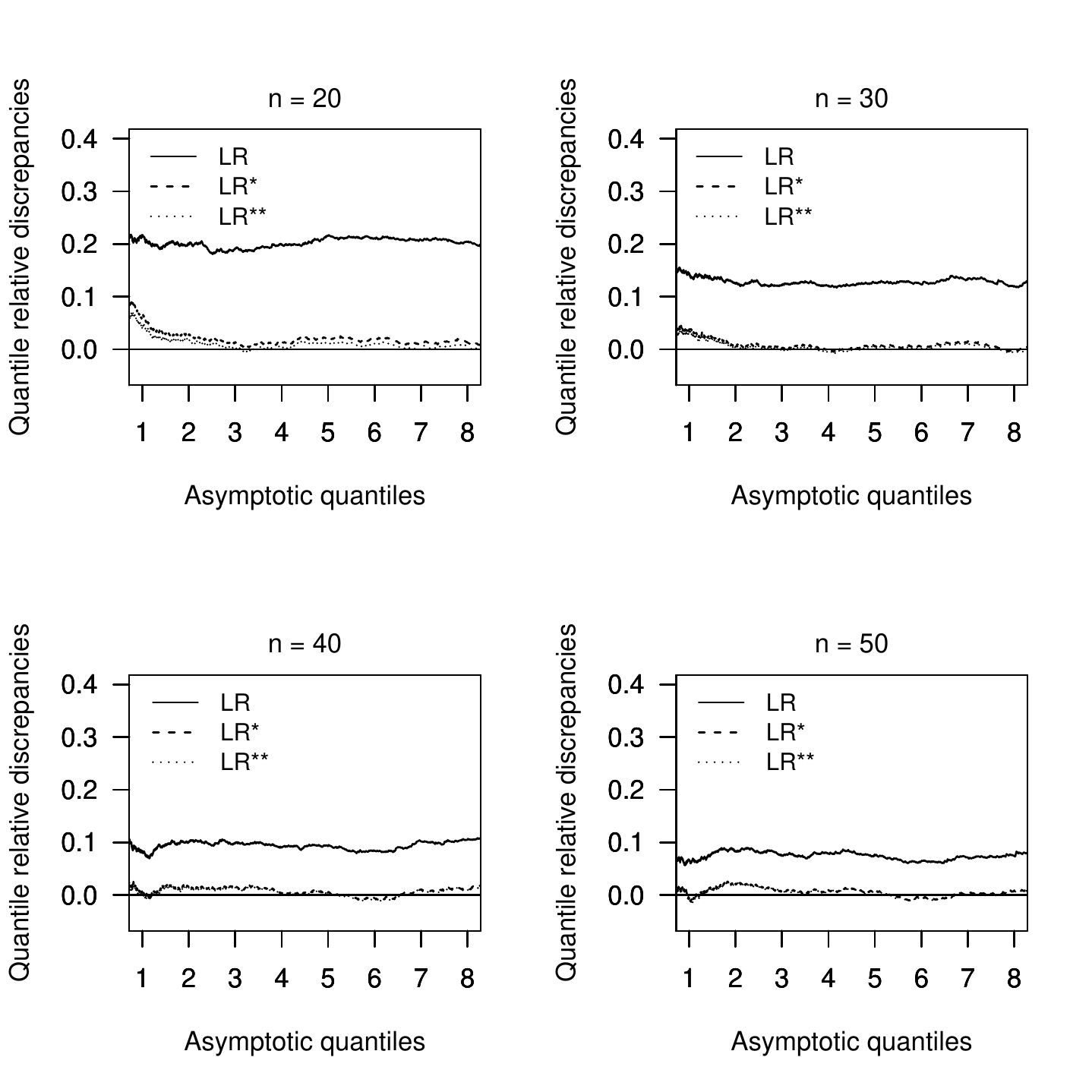}
\caption{Quantile relative discrepancies: $p = 3$, $q = 3$, and power exponential distribution for ${\bm Z}_i$.}
\label{EP3q3}
\end{center}
\end{figure}

Figures \ref{tStudp4q3}, \ref{Normp4q4} and \ref{EP3q3} depict curves of quantile relative discrepancies \emph{versus} the correspondent asymptotic quantiles for the statistics $LR$, $LR_a^*$, and $LR_a^{**}$ considering different sample sizes and $p = 4$ and $q = 3$ under a Student-$t$ distribution (Figure \ref{tStudp4q3}), $p = q = 4$ under a normal distribution (Figure \ref{Normp4q4}), and $p = q = 3$  under a power exponential distribution (Figure \ref{EP3q3}). The closer to zero a curve is, the better is the asymptotic approximation. As expected, the three figures show that the curves of relative quantile discrepancies for the statistics  $LR_a^{*}$ and $LR_a^{**}$ are closer to the horizontal axis of the zero ordinate than the curve based on $LR$.

Tables 4 and 5 present rejection rates obtained under alternative hypotheses for $n=20$ and $n=40$, respectively. Recall that the null hypothesis is 
${\mathcal H}_{0}: {\bm \psi} = {\bm 0}$, and that $\psi=(\beta_1,\ldots,\beta_q)$. Here, the rejection rates were obtained assuming that
${\bm \psi} = (\eta, \ldots,\eta)$, with  ${\eta}$ varying from $0.1$ to $1.5$, and $p = q = 2$. The original likelihood ratio test
was not included in this study, since our simulations indicated that it is oversized. The figures in Tables 4 and 5 indicate that 
the modified tests present similar power, the power of the test that uses $LR^*_a$ being equal to or slightly greater than that of the power of the test based on  $LR^{**}_a$.

We conclude that the modified versions of the likelihood ratio test perform better than the original test for small and moderate sample sizes.  Although the tests based on the modified statistics, $LR_a^*$ and $LR_a^{**}$, present similar results, the test based on $LR_a^{**}$ has a slightly better performance (in the majority of the cases) than the one based on $LR_a^*$. 

\begin{table}[!htp]
\hspace{1.6cm} {\caption{Rejection rates of  ${\mathcal H}_0$ when $q = 2$ and $p = 2, 3, 4$; normal distribution.}} 
\vspace{0.5cm} \centering
\begin{tabular}{cccccccccccc} 
\hline
\multicolumn{12}{c}{$p = 2$}\\
\hline
&\multicolumn{3}{c}{$\gamma = 1\%$}&&\multicolumn{3}{c}{$\gamma = 5\%$}&&\multicolumn{3}{c}{$\gamma = 10\%$}
\\\cline{2-4}\cline{6-8}\cline{10-12}
{$n$}& $LR$ & $LR_a^*$& $LR_a^{**}$ && $LR$  & $LR_a^*$& $LR_a^{**}$ && $LR$  & $LR_a^*$& $LR_a^{**}$ \\\hline 
  20 &2.1&1.1&1.0&&8.2&5.4&5.2&&14.7&10.6&10.4\\ 
  30 &1.7&1.1&1.1&&6.9&5.2&5.1&&13.0&10.6&10.5\\ 
  40 &1.5&1.1&1.1&&6.2&4.9&4.9&&11.9&10.0&10.0\\   
  50 &1.5&1.1&1.1&&6.2&5.2&5.2&&11.7&10.3&10.3\\ 
\hline
\multicolumn{12}{c}{$p = 3$}\\
\hline
&\multicolumn{3}{c}{$\gamma = 1\%$}&&\multicolumn{3}{c}{$\gamma = 5\%$}&&\multicolumn{3}{c}{$\gamma = 10\%$}
\\\cline{2-4}\cline{6-8}\cline{10-12}
{$n$}& $LR$ & $LR_a^*$& $LR_a^{**}$ && $LR$  & $LR_a^*$& $LR_a^{**}$ && $LR$  & $LR_a^*$& $LR_a^{**}$ \\\hline 
  20 &2.6&1.1&1.0&&9.3&5.6&5.3&&16.0&10.5&10.2\\ 
  30 &2.2&1.2&1.2&&7.8&5.3&5.1&&13.9&10.6&10.4\\ 
  40 &1.5&1.0&1.0&&6.9&5.1&5.0&&12.8&10.4&10.3\\ 
  50 &1.6&1.1&1.1&&6.8&5.4&5.4&&12.1&10.4&10.4\\ 
\hline
\multicolumn{12}{c}{$p = 4$}\\
\hline
&\multicolumn{3}{c}{$\gamma = 1\%$}&&\multicolumn{3}{c}{$\gamma = 5\%$}&&\multicolumn{3}{c}{$\gamma = 10\%$}
\\\cline{2-4}\cline{6-8}\cline{10-12}
{$n$}& $LR$ & $LR_a^*$& $LR_a^{**}$ && $LR$  & $LR_a^*$& $LR_a^{**}$ && $LR$  & $LR_a^*$& $LR_a^{**}$ \\\hline 
20 &3.4&1.2&1.1&&11.3&6.1&5.7&&18.8&11.7&11.2 \\
30 &2.1&1.1&1.0&& 8.4&5.4&5.2&&15.0&10.6&10.4 \\
40 &1.6&1.1&1.0&& 7.3&5.1&5.1&&13.9&10.7&10.5 \\
50 &1.5&0.9&0.9&& 6.5&4.9&4.8&&12.7&10.3&10.3 \\
\hline
\end{tabular}
\label{tab.1}
\end{table}

\begin{table}[!htp]
\hspace{1.6cm} {\caption{Rejection rates of ${\mathcal H}_0$ when $q = 2$ and $p = 2, 3, 4$; Student-$t$ distribution ($\nu = 5$).}} 
\vspace{0.5cm} \centering
\begin{tabular}{cccccccccccc} 
\hline
\multicolumn{12}{c}{$p = 2$}\\
\hline
&\multicolumn{3}{c}{$\gamma = 1\%$}&&\multicolumn{3}{c}{$\gamma = 5\%$}&&\multicolumn{3}{c}{$\gamma = 10\%$}
\\\cline{2-4}\cline{6-8}\cline{10-12}
{$n$}& $LR$ & $LR_a^*$& $LR_a^{**}$ && $LR$  & $LR_a^*$& $LR_a^{**}$ && $LR$  & $LR_a^*$& $LR_a^{**}$ \\\hline 
  20 &1.7&0.9&0.8&&7.4&4.8&4.6&&13.6&10.1& 9.8\\ 
  30 &1.4&0.9&0.9&&6.5&4.7&4.6&&12.1& 9.9& 9.8\\  
  40 &1.4&1.0&1.0&&6.1&4.9&4.9&&11.5& 9.7& 9.6\\  
  50 &1.3&1.0&1.0&&6.2&5.1&5.1&&11.5&10.2&10.2\\  
\hline
\multicolumn{12}{c}{$p = 3$}\\
\hline
&\multicolumn{3}{c}{$\gamma = 1\%$}&&\multicolumn{3}{c}{$\gamma = 5\%$}&&\multicolumn{3}{c}{$\gamma = 10\%$}
\\\cline{2-4}\cline{6-8}\cline{10-12}
{$n$}& $LR$ & $LR_a^*$& $LR_a^{**}$ && $LR$  & $LR_a^*$& $LR_a^{**}$ && $LR$  & $LR_a^*$& $LR_a^{**}$ \\\hline 
  20 &2.5&1.2&1.1&&9.0&5.5&5.3&&15.7&10.9&10.5 \\ 
  30 &2.1&1.2&1.1&&7.4&5.4&5.3&&13.5&10.4&10.3\\  
  40 &1.5&1.0&0.9&&6.3&4.8&4.8&&12.0& 9.8& 9.7\\  
  50 &1.4&1.0&1.0&&6.1&4.8&4.8&&11.8&10.1&10.0\\  
\hline
\multicolumn{12}{c}{$p = 4$}\\
\hline
&\multicolumn{3}{c}{$\gamma = 1\%$}&&\multicolumn{3}{c}{$\gamma = 5\%$}&&\multicolumn{3}{c}{$\gamma = 10\%$}
\\\cline{2-4}\cline{6-8}\cline{10-12}
{$n$}& $LR$ & $LR_a^*$& $LR_a^{**}$ && $LR$  & $LR_a^*$& $LR_a^{**}$ && $LR$  & $LR_a^*$& $LR_a^{**}$ \\\hline 
  20 &2.6&1.3&1.0&&9.2&5.6&5.1&&16.1&10.9&10.2\\ 
  30 &2.0&0.9&0.9&&7.5&5.2&5.1&&13.9&10.2&10.0\\  
  40 &1.9&1.2&1.1&&7.2&5.5&5.5&&13.1&10.5&10.4\\  
  50 &1.5&1.0&1.0&&6.6&5.1&5.1&&12.3&10.2&10.1\\  
\hline
\end{tabular}
\label{tab.2}
\end{table}

\begin{table}[!htp]
\hspace{1.6cm} {\caption{Rejection rates of ${\mathcal H}_0$ when $q = 2$ and $p = 2, 3, 4$; power exponential distribution ($\lambda = 0.6$).}} 
\vspace{0.5cm} \centering
\begin{tabular}{cccccccccccc} 
\hline
\multicolumn{12}{c}{$p = 2$}\\
\hline
&\multicolumn{3}{c}{$\gamma = 1\%$}&&\multicolumn{3}{c}{$\gamma = 5\%$}&&\multicolumn{3}{c}{$\gamma = 10\%$}
\\\cline{2-4}\cline{6-8}\cline{10-12}
{$n$}& $LR$ & $LR_a^*$& $LR_a^{**}$ && $LR$  & $LR_a^*$& $LR_a^{**}$ && $LR$  & $LR_a^*$& $LR_a^{**}$ \\\hline 
  20 &1.9&1.0&1.0&&7.2&5.0&4.9&&13.5& 9.8& 9.6\\ 
  30 &1.4&1.0&0.9&&6.4&4.7&4.7&&12.2& 9.9& 9.9\\  
  40 &1.3&1.0&1.0&&5.9&4.9&4.9&&11.5&10.0& 9.9\\  
  50 &1.4&1.2&1.2&&6.2&5.1&5.1&&11.8&10.4&10.4\\  
\hline
\multicolumn{12}{c}{$p = 3$}\\
\hline
&\multicolumn{3}{c}{$\gamma = 1\%$}&&\multicolumn{3}{c}{$\gamma = 5\%$}&&\multicolumn{3}{c}{$\gamma = 10\%$}
\\\cline{2-4}\cline{6-8}\cline{10-12}
{$n$}& $LR$ & $LR_a^*$& $LR_a^{**}$ && $LR$  & $LR_a^*$& $LR_a^{**}$ && $LR$  & $LR_a^*$& $LR_a^{**}$ \\\hline 
  20 &2.4&1.2&1.1&&9.3&5.5&5.3&&15.5&10.9&10.6\\ 
  30 &1.9&1.1&1.1&&7.5&5.4&5.3&&13.1&10.1& 9.9\\  
  40 &1.8&1.1&1.1&&6.7&4.9&4.9&&12.4&10.3&10.2\\  
  50 &1.5&1.2&1.2&&6.2&4.9&4.9&&12.5&10.3&10.3\\  
\hline
\multicolumn{12}{c}{$p = 4$}\\
\hline
&\multicolumn{3}{c}{$\gamma = 1\%$}&&\multicolumn{3}{c}{$\gamma = 5\%$}&&\multicolumn{3}{c}{$\gamma = 10\%$}
\\\cline{2-4}\cline{6-8}\cline{10-12}
{$n$}& $LR$ & $LR_a^*$& $LR_a^{**}$ && $LR$  & $LR_a^*$& $LR_a^{**}$ && $LR$  & $LR_a^*$& $LR_a^{**}$ \\\hline 
  20 &2.9&1.2&1.0&&10.3&5.9&5.4&&17.3&11.6&10.9\\ 
  30 &2.1&1.2&1.2&& 7.8&5.1&4.9&&14.0&10.2&10.0\\  
  40 &1.9&1.1&1.1&& 6.9&4.8&4.8&&12.8&10.0& 9.8\\  
  50 &1.5&1.1&1.0&& 6.6&5.1&5.1&&11.9& 9.8& 9.7\\  
\hline
\end{tabular}
\label{tab.3}
\end{table}

\begin{table}[h]
\hspace{1.6cm} {\caption{\scriptsize Power of the modified tests when $n=20$ and $p = q = 2$; normal, Student-$t$ ($\nu=5$) and power exponential ($\lambda=0.6$) distributions.}} 
\vspace{0.5cm} \centering
\begin{tabular}{cccccc|ccccc|ccccc}
\hline 
\vspace{0.1cm}
&\multicolumn{5}{c}{normal}&\multicolumn{5}{c}{Student-$t$}&\multicolumn{5}{c}{power exponential}
\\
\cline{2-6}\cline{7-11}\cline{12-16}
&\multicolumn{2}{c}{$\gamma = 5\%$}&&\multicolumn{2}{c}{$\gamma = 10\%$} &\multicolumn{2}{|c}{$\gamma = 5\%$}&&\multicolumn{2}{c|}{$\gamma = 10\%$} &\multicolumn{2}{c}{$\gamma = 5\%$}&&\multicolumn{2}{c}{$\gamma = 10\%$}
\\\cline{2-3}\cline{5-6} \cline{7-8}\cline{10-11} \cline{12-13}\cline{15-16}
{${\bm \eta}$}&$LR_a^*$&$LR_a^{**}$&&$LR_a^*$&$LR_a^{**}$    &$LR_a^*$&$LR_a^{**}$&&$LR_a^*$&$LR_a^{**}$    &$LR_a^*$&$LR_a^{**}$&&$LR_a^*$&$LR_a^{**}$\\\hline
  0.1    &5.4  &5.2  &&10.6 &10.4 &4.8  &4.6  &&10.1 &9.8  &5.0  &4.9  &&9.8  &9.6  \\  
  0.2    &6.2  &6.1  &&12.0 &11.8 &5.6  &5.5  &&11.4 &11.2 &5.6  &5.6  &&11.3 &11.1  \\  	   
  0.3    &9.5  &9.3  &&16.4 &16.2 &8.1  &8.0  &&14.8 &14.5 &8.3  &8.2  &&15.1 &14.8  \\       
  0.4    &14.0 &13.8 &&23.3 &23.0 &12.3 &12.1 &&20.8 &20.6 &12.8 &12.6 &&21.3 &21.0  \\       
  0.5    &22.3 &22.0 &&33.4 &33.1 &18.4 &18.1 &&28.7 &28.5 &19.3 &19.1 &&29.9 &29.6  \\  
  0.6    &30.9 &30.6 &&43.9 &43.5 &26.9 &26.6 &&38.6 &38.3 &27.4 &27.1 &&39.2 &38.9  \\       
  0.7    &43.0 &42.6 &&56.5 &56.2 &35.3 &34.9 &&48.6 &48.2 &37.0 &36.5 &&50.0 &49.7  \\  		   
  0.8    &53.9 &53.6 &&66.8 &66.4 &44.7 &44.4 &&58.1 &57.8 &46.9 &46.5 &&60.5 &60.2  \\  		   
  0.9    &64.9 &64.5 &&76.0 &75.8 &57.3 &57.0 &&69.8 &69.6 &56.5 &56.1 &&69.3 &69.1  \\  		   
  1.0    &72.8 &72.5 &&82.4 &82.2 &64.9 &64.5 &&76.3 &76.0 &67.5 &67.2 &&78.4 &78.1  \\  		   
  1.1    &81.5 &81.1 &&88.4 &88.2 &71.1 &70.7 &&81.0 &80.8 &74.5 &74.2 &&84.0 &83.7  \\			
  1.2    &86.8 &86.7 &&92.4 &92.2 &79.0 &78.8 &&87.2 &87.1 &81.1 &80.9 &&89.4 &89.2  \\
  1.3    &90.9 &90.8 &&95.4 &95.3 &84.4 &84.2 &&91.0 &90.8 &86.0 &85.7 &&91.9 &91.7  \\	 
  1.4    &93.5 &93.4 &&96.8 &96.8 &90.1 &89.9 &&94.9 &94.9 &90.5 &90.3 &&95.4 &95.3  \\	 
  1.5    &97.4 &97.3 &&99.6 &99.5 &93.2 &93.0 &&96.9 &96.9 &92.3 &92.2 &&96.4 &96.3  \\
\hline
\end{tabular}
\label{tab.4}
\end{table}

\begin{table}[h]
\hspace{1.6cm} {\caption{\scriptsize Power of the modified tests when $n=40$ and $p = q = 2$; normal, Student-$t$ ($\nu=5$) and power exponential ($\lambda=0.6$) distributions.}} 
\vspace{0.5cm} \centering
\begin{tabular}{cccccc|ccccc|ccccc}
\hline 
\vspace{0.1cm}
&\multicolumn{5}{c}{normal}&\multicolumn{5}{c}{Student-$t$}&\multicolumn{5}{c}{power exponential}
\\
\cline{2-6}\cline{7-11}\cline{12-16}
&\multicolumn{2}{c}{$\gamma = 5\%$}&&\multicolumn{2}{c}{$\gamma = 10\%$} &\multicolumn{2}{|c}{$\gamma = 5\%$}&&\multicolumn{2}{c|}{$\gamma = 10\%$} &\multicolumn{2}{c}{$\gamma = 5\%$}&&\multicolumn{2}{c}{$\gamma = 10\%$}
\\\cline{2-3}\cline{5-6} \cline{7-8}\cline{10-11} \cline{12-13}\cline{15-16}
{${\bm \eta}$}&$LR_a^*$&$LR_a^{**}$&&$LR_a^*$&$LR_a^{**}$    &$LR_a^*$&$LR_a^{**}$&&$LR_a^*$&$LR_a^{**}$    &$LR_a^*$&$LR_a^{**}$&&$LR_a^*$&$LR_a^{**}$\\\hline
  0.1    &7.4	 &7.4	 &&13.6  &13.5  &6.5  &6.4  &&12.4	 &12.4  &6.7	 &6.6    &&12.6	 &12.6  \\  
  0.2    &14.1	 &14.0	 &&23.3  &23.3  &12.2 &12.1 &&20.7	 &20.6  &12.6	 &12.5   &&20.8	 &20.7   \\  	    
  0.3    &26.5	 &26.4	 &&38.2  &38.2  &22.5 &22.4 &&33.1	 &33.0  &23.0	 &22.9   &&33.7	 &33.7   \\       
  0.4    &43.0	 &42.9	 &&55.9  &55.9  &35.8 &35.8 &&48.5	 &48.4  &36.9	 &36.8   &&49.7	 &49.7   \\       
  0.5    &60.7	 &60.6	 &&72.4  &72.3  &51.4 &51.4 &&64.4	 &64.3  &53.6	 &53.5   &&65.9	 &65.8   \\  
  0.6    &76.1	 &76.0	 &&85.2  &85.1  &66.8 &66.7 &&77.3	 &77.3  &68.8	 &68.7   &&78.4	 &78.3   \\       
  0.7    &87.5	 &87.4	 &&92.9  &92.9  &79.1 &79.1 &&86.9	 &86.9  &80.3	 &80.2   &&87.8	 &87.8   \\  		   
  0.8    &94.0	 &93.9	 &&97.1  &97.1  &88.1 &88.1 &&93.2	 &93.2  &88.9	 &88.9   &&93.7	 &93.7   \\  		   
  0.9    &97.6	 &97.6	 &&98.8  &98.8  &93.7 &93.7 &&96.8	 &96.8  &94.0	 &94.0   &&97.1	 &97.0   \\  		   
  1.0    &99.0	 &99.0	 &&99.5  &99.5  &96.9 &96.9 &&98.5	 &98.5  &97.2	 &97.1   &&98.7	 &98.7   \\  		   
  1.1    &99.6	 &99.6	 &&99.9  &99.9  &98.5 &98.5 &&99.3	 &99.3  &98.8	 &98.8   &&99.5	 &99.5   \\			
  1.2    &99.9	 &99.9	 &&100.0 &100.0 &99.3 &99.3 &&99.7	 &99.7  &99.5	 &99.5   &&99.8	 &99.8   \\
  1.3    &100.0  &100.0  &&100.0 &100.0 &99.7 &99.7 &&99.9	 &99.9  &99.8	 &99.8   &&99.9	 &99.9   \\	 
  1.4    &100.0  &100.0  &&100.0 &100.0 &99.9 &99.9 &&99.9	 &99.9  &100.0   &100.0  &&100.0 &100.0  \\	 
  1.5    &100.0  &100.0  &&100.0 &100.0 &99.9 &99.9 &&100.0  &100.0 &100.0   &100.0  &&100.0 &100.0  \\
\hline
\end{tabular}
\label{tab.5}
\end{table}

\section{Application}
\label{sec5}

The WHO MONICA Project was established in the early 1980s, and has been monitoring trends in cardiovascular diseases since then. One of the goals of this project is to relate changes in known risk factors ($x$) with the trends in cardiovascular mortality and coronary heart disease ($y$). In this paper, we consider a data set from the WHO MONICA Project that was considered in Kulathinal et al. \cite{Kulathinaletal2002}, Patriota et al. \cite{Patriotaetal2009}, and Cao et al. \cite{Caoetal2012}. This data set was first analyzed under normal distributions for the marginals of the random errors (\cite{Kulathinaletal2002}, \cite{Patriotaetal2009}), and thereafter under scale-mixture of normal distributions for the marginals of the random errors (Cao et al. \cite{Caoetal2012}). The approach used in the present paper is different from the latter, since here we consider an elliptical joint distribution for the vector of the random errors. Although the class of scale-mixture of normal distributions is a special case of the the elliptical distributions, our proposal does not extend to the class proposed by Cao et al. \cite{Caoetal2012}. This is because they assumed that the {errors' distributions} are independent, while we assume that they are uncorrelated but not independent. For our proposal, the errors will only be independent under normality.

The data set considered here corresponds to the data collected for men ($n=36$). As describe in Kulathinal et al. \cite{Kulathinaletal2002}, the data are trends of the annual change in event rate $(y)$ and trends of the risk scores ($x$). The risk score is defined as a linear combination of smoking status, systolic blood pressure, body mass index, and total cholesterol level. A follow-up study using proportional hazards models was employed to derive its coefficients, and provides the observed risk score and its estimated variance. Therefore, the observed response variable, $Y$, is the average annual change in event rate (\%) and the observed covariate, $X$, is the observed risk score (\%). For these data, we use model (\ref{E.4.1}) with $m=p=1$ and zero covariance between the errors $e$ and $u$.

Table \ref{Res.men} gives the two corrected statistics and the original likelihood ratio statistic for testing ${\cal H}_0: \beta = 1$ ($p$-values are given in parentheses). We considered the full sample ($n=36$) and a randomly chosen sub-sample of $n=20$ observations. We notice that for the full sample all the tests have similar $p$-values, and reject the null hypothesis at the 5\% nominal level. 
For the sub-sample, however, the null hypothesis is rejected at the 10\% nominal level when one uses either of the corrected statistics, but ${\cal H}_0$ is not rejected when the usual likelihood ratio test is employed.

\begin{table}[!htp]\caption{Corrected and uncorrected likelihood ratio statistics for testing ${\cal H}_0:\beta = 1$; $p$-values are given in parentheses. }\label{Res.men}
\[
\begin{array}{ccccc}\hline
          & {\mbox{Full sample}}   & & {\mbox{Sub-sample}}\\
\hline
LR  	    &4.5037 \ (0.0338)      & &	2.4160 \ (0.1201)\\
LR_a^* 	  &4.5701 \ (0.0325)      & &	3.6583 \ (0.0558)\\
LR_a^{**}	&4.5698 \ (0.0325)      & &	3.5299 \ (0.0603)\\
\hline
\end{array}
\]
\end{table}

\section{Concluding remarks}
\label{sec6}

We studied the issue of testing hypotheses in heteroskedastic multivariate regression models with measurement error under the elliptical class of distributions.  For this class of distributions, the exact distribution of the likelihood ratio statistic under the null hypothesis is not easily attained and the usual asymptotic approximation to the chi-square distribution is required. However, for small sample sizes, this approximation may be very inaccurate and modified versions proposed by Skovgaard \cite{Skovgaard2001} can be used to obtain more accurate approximations under the null hypothesis. The main difficulty of this approach is in finding a suitable ancillary statistic. In this work, we find this ancillary statistic and all the quantities required to compute Skovgaard's adjustment. Simulation results indicate that the adjusted tests are more reliable than the original likelihood ratio test for small sample sizes. A real data application illustrates the usefulness of the results presented in this paper.

\section*{Acknowledgments}

We gratefully acknowledge the financial the support from FAPESP and CNPq. The authors thank Dr.~Kari Kuulasmaa (National Public Health
Institute, Finland) for kindly supplying the data, and the anonymous referee for helpful comments and suggestions.

\begin{appendix}
{\small

\section*{Appendix A. Observed information matrix}
\label{secA}

The first-order derivative of  (\ref{E.4.8}) with respect to $\theta_{j}$, with $j = 1, 2, \ldots, s$, is
\begin{eqnarray}\label{C.1.2}
\frac{\partial \ell(\bm{\theta})}{\partial \theta_{j}} = \sum_{i=1}^{n} \left[-\frac{1}{2} \tr\left({\Omega}_i^{-1} \Omega_{{\theta_j}}\right) - W_{p_0}({\bm d}_i^\top \Omega_i^{-1} {\bm d}_i) \left({\bm d}_i^\top {\Omega}_i^{-1} {\Omega}_{\theta_j} {\Omega}_i^{-1} {\bm d}_i + 2 \bm{\mu}_{\theta_j} {\Omega}_i^{-1} {\bm d}_i\right)\right],
\end{eqnarray}
where $W_{p_0}(u) = {\partial\log p_0(u)}/{\partial u}$ .
Define
\[
F_{\beta}^{(j)} = \frac{\partial \bm{\beta}}{\partial \theta_j}, \quad F_{\alpha}^{(j)} = \frac{\partial \bm{\alpha}}{\partial \theta_j}, \quad F_{\mu_x}^{(j)} = \frac{\partial \bm{\mu}_x}{\partial \theta_j},  
\]
\[
F_{\Sigma_q}^{(j)} = \frac{\partial \Sigma_q}{\partial \theta_j}, \quad \mbox{and} \quad F_{\Sigma_x}^{(j)} = \frac{\partial \Sigma_x }{\partial \theta_j}.  
\]
Then, the vector $\bm{\mu}_{\theta_j} = \partial {\bm{\mu}}/\partial\theta_j$ and matrix  $\Omega_{\theta_j} = \partial\Omega_i/\partial\theta_j$ are given, respectively, by 
\begin{eqnarray*}\label{C.1.3}
\bm{\mu}_{\theta_j} = \begin{pmatrix} F_{\beta}^{(j)}\mu_{x}\\ {\bm 0} \end{pmatrix} 
                    + \begin{pmatrix} F_{\alpha}^{(j)}\\ {\bm 0} \end{pmatrix} 
                    + \begin{pmatrix} \bm{\beta}F_{\mu_x}^{(j)}\\ F_{\mu_x}^{(j)} \end{pmatrix}
\end{eqnarray*}
and
\begin{eqnarray*}\label{C.1.4}
\Omega_{\theta_j} =  c^{-1} \begin{pmatrix} 
                     F_{\beta}^{(j)}\Sigma_x \bm{\beta}^\top + \bm{\beta} \Sigma_x F_{\beta}^{\top(j)} & F_{\beta}^{(j)}\Sigma_x\\
                     \Sigma_x F_{\beta}^{\top(j)}                                                      & {\bm 0} 
                     \end{pmatrix}
                     + 
                     c^{-1} \begin{pmatrix} 
                     F_{\Sigma_q}^{(j)} & {\bm 0}\\ 
                     {\bm 0} & {\bm 0}
                     \end{pmatrix}
                     + 
                     c^{-1} \begin{pmatrix} 
                     \bm{\beta}F_{\Sigma_x}^{(j)}\bm{\beta}^{(j)} & \bm{\beta}F_{\Sigma_x}^{(j)}\\ 
                     F_{\Sigma_x}^{(j)}\bm{\beta}^\top & F_{\Sigma_x}^{(j)}
                     \end{pmatrix},
\end{eqnarray*}
for $j, k = 1, 2, \ldots, s$. Note that the derivative presented in (\ref{C.1.2}) forms the score vector. 
The element $(j,k)$ of the observed information matrix $J = J(\bm{\theta})$ is $J_{\theta_j\theta_k} = -{\partial^2 \ell(\bm{\theta})}/{\partial \theta_{j} \partial \theta_{k}}$, i.e.,
\begin{equation*}\label{C.1.10}
\begin{split}
J_{\theta_j\theta_k} &=  \sum_{i=1}^{n}\Bigg\{-\frac{1}{2} \tr({\Omega}_i^{-1} {\Omega}_{\theta_k} {\Omega}_i^{-1} \Omega_{\theta_j}) + \frac{1}{2} \tr({\Omega}_i^{-1} \Omega_{\theta_j \theta_k}) - W_{p_0}'({\bm d}_i^\top \Omega_i^{-1} {\bm d}_i)\Big[{\bm d}_i^\top {\Omega}_i^{-1} {\Omega}_{\theta_k} {\Omega}_i^{-1} {\bm d}_i 
\\
&+ 2 \bm{\mu}_{\theta_k}^\top  {\Omega}_i^{-1} {\bm d}_i \Big] \Big[{\bm d}_i^\top {\Omega}_i^{-1} {\Omega}_{\theta_j} {\Omega}_i^{-1} {\bm d}_i + 2\bm{\mu}_{\theta_j}^\top  {\Omega}_i^{-1} {\bm d}_i \Big] + W_{p_0}({\bm d}_i^\top \Omega_i^{-1} {\bm d}_i) \Big[{\bm d}_i^\top \Big( 2{\Omega}_i^{-1} {\Omega}_{\theta_k} {\Omega}_i^{-1} {\Omega}_{\theta_j} \Omega_i^{-1} 
\\
&- \Omega_i^{-1} \Omega_{\theta_j\theta_k}\Omega_i^{-1}\Big) {\bm d}_i + 2\bm{\mu}_{\theta_k}^\top {\Omega}_i^{-1} {\Omega}_{\theta_j} {\Omega}_i^{-1} {\bm d}_i + 2\bm{\mu}_{\theta_j}^\top {\Omega}_i^{-1} {\Omega}_{\theta_k} {\Omega}_i^{-1}  {\bm d}_i - 2\bm{\mu}_{\theta_j\theta_k}^\top  {\Omega}_i^{-1} {\bm d}_i
 + 2 \bm{\mu}_{\theta_j}^\top {\Omega}_i^{-1} \bm{\mu}_{\theta_k}\Big]\Bigg\},
\end{split}
\end{equation*}
for $j, k = 1, 2, \ldots, s$, where $W_{p_0}'(u) = {\partial W_{p_0}(u)}/{\partial u}$. The vector $\bm{\mu}_{\theta_j\theta_k} = {\partial\bm{\mu}_{\theta_j}}/{\partial \theta_k}$ is 
\begin{eqnarray*}\label{C.1.7}
\bm{\mu}_{\theta_j\theta_k} = \bm{\mu}_{\theta_k\theta_j} = 
       \begin{pmatrix} 
	F_{\beta}^{(j)}F_{\mu_{x}}^{(k)}\\ 
	{\bm 0} 
	\end{pmatrix}+
 \begin{pmatrix} 
	F_{\beta}^{(k)}F_{\mu_{x}}^{(j)}\\ 
	{\bm 0} 
	\end{pmatrix},  
\end{eqnarray*} for $j,k = 1,2, \ldots, s$. The matrix $\Omega_{\theta_j\theta_k} = {\partial\Omega_{\theta_j}}/{\partial \theta_k}$ is symmetric and has dimension $s\times s$, with $j, k = 1, 2, \ldots, s$, and is given by
\begin{eqnarray*}\label{C.1.8}
\Omega_{\theta_j\theta_k} =c^{-1}
         \begin{pmatrix} 
         F_{\beta}^{(j)}F_{\Sigma_x}^{(k)} \bm{\beta}^\top + \bm{\beta} F_{\Sigma_x}^{(k)} F_{\beta}^{\top(j)} & F_{\beta}^{(j)}F_{\Sigma_x}^{(k)}\\
         F_{\Sigma_x}^{(k)} F_{\beta}^{\top(j)}                                                      & {\bm 0} 
         \end{pmatrix}+c^{-1}
 \begin{pmatrix} 
         F_{\beta}^{(k)}F_{\Sigma_x}^{(j)} \bm{\beta}^\top + \bm{\beta} F_{\Sigma_x}^{(j)} F_{\beta}^{\top(k)} & F_{\beta}^{(k)}F_{\Sigma_x}^{(j)}\\
         F_{\Sigma_x}^{(j)} F_{\beta}^{\top(k)}                                                      & {\bm 0} 
         \end{pmatrix}, \\
\end{eqnarray*} for $j,k = 1, 2, \ldots, s$. 

In matrix notation, the observed information matrix is given by
\begin{equation*}\label{C.1.11}
J = \frac{1}{2} T - R^\top M - V^\top Q,
\end{equation*}
\newpage
where the expressions for the elements of  $T$, $R$, $M$, $V$, and  $Q$ are presented in Section \ref{sec3}.

\vspace{-0.5cm}
\section*{Appendix B. Derivatives with respect to the data}
\label{secB}
We present the derivatives with respect to the sample space $\ell'$, $U'$, and ${\:{\J}}$. Consider the ancillary statistic ${\bm a} = ({\bm a}_1^\top, {\bm a}_2^\top$, $\ldots, {\bm a}_n^\top)^\top$ defined in  $(\ref{E.20})$, where ${\bm a}_i = {\widehat P}_i^{-1}\left({\bm z}_i - {\bm{\widehat \mu}}\right)$. Inserting ${\bm z}_i = {\widehat P}_i {\bm a}_i + \bm{\widehat \mu}$ in the log-likelihood function $(\ref{E.4.8})$, we obtain 
\begin{eqnarray*}\label{C.2.1}
\ell(\bm{\theta};\bm{\widehat \theta},{\bm a}) = \sum_{i=1}^{n} \left\{-\frac{1}{2} \log |\Omega_i| + \log p_0\left[({\widehat P}_i {\bm a}_i + \bm{\widehat \mu} - \bm{\mu})^\top \Omega_i^{-1} ({\widehat P}_i {\bm a}_i + \bm{\widehat \mu} - \bm{\mu})\right]\right\}.
\end{eqnarray*}
Therefore,
the $j$-th element of $\ell' = {\partial \ell(\bm{\theta};\bm{\widehat \theta},a)}/{\partial \bm{\widehat \theta}}$ is
\begin{equation*}\label{C.2.5}
{\ell_j'} = 2 \sum_{i=1}^{n} W_{p_0}\left({\bm d}_i^\top \Omega_i^{-1} {\bm d}_i\right) \left({\widehat P}_{\theta_j} {\bm a}_i + \bm{\widehat \mu}_{\theta_j}\right)^\top \Omega_i^{-1} {\bm d}_i.
\end{equation*}
Furthermore, the $(j,k)$-th element of $U'= {\partial^2 \ell(\bm{\theta};\bm{\widehat \theta},{\bm a})}/{\partial \bm{\widehat \theta} \partial{\bm{\theta}^\top}}$ is 
\begin{equation*}\label{C.2.6}
\begin{split}
U_{jk}'&= 2\sum_{i=1}^{n}\Bigg\{W_{p_0}\left({\bm d}_i^\top \Omega_i^{-1}{\bm d}_i\right)\left[-\left({\widehat P}_{\theta_k} {\bm a}_i + \bm{\widehat \mu}_{\theta_k}\right)^\top {\Omega}_i^{-1} {\Omega}_{\theta_j} {\Omega}_i^{-1} {\bm d}_i 
 - {\bm{\mu}}_{\theta_j}^\top \Omega_i^{-1} \left({\widehat P}_{\theta_k} {\bm a}_i 
+ \bm{\widehat \mu}_{\theta_k}\right)\right]
\\
& + W'_{p_0}\left({\bm d}_i^\top \Omega_i^{-1} {\bm d}_i\right)
\left({\widehat P}_{\theta_k} {\bm a}_i + \bm{\widehat \mu}_{\theta_k}\right)^\top \Omega_i^{-1} {\bm d}_i\left(-{\bm d}_i^\top {\Omega}_i^{-1} {\Omega}_{\theta_j} {\Omega}_i^{-1}{\bm d}_i 
 - 2 \bm{\mu}_{\theta_k}^\top {\Omega}_i^{-1}{\bm d}_i\right)\Bigg\},
\end{split}
\end{equation*}
where ${\widehat P}_{\theta_j} = {\partial {\widehat P}_i}/{{\widehat \theta_j}}$, for $j, k = 1, 2, \ldots, s$. We also have that the $(j,k)$-th element of ${\:\:{\J}}$
is 
\begin{equation*}\label{C.2.7}
\begin{split}
{\:\:{\J}}_{jk}&= 2\sum_{i=1}^{n}\Bigg\{W_{p_0}\left({\bm{\widehat d}}_i^\top {\widehat \Omega}_i^{-1} {\bm{\widehat d}}_i\right)\left[-\left({\widetilde P}_{\theta_k} {\bm a}_i +  \bm{\widetilde \mu}_{\theta_k}\right)^\top {{\widetilde\Omega}}_i^{-1} {{\widetilde\Omega}}_{\theta_j} {{\widetilde\Omega}}_i^{-1}{\widetilde P}_i {\bm a}_i - {\bm{\widetilde\mu}}_{\theta_j}^\top {\widetilde\Omega}_i^{-1} \left({\widetilde P}_{\theta_k} {\bm a}_i + \bm{\widetilde \mu}_{\theta_k}\right)\right]
\\
&+ W'_{p_0}\left({\bm{\widehat d}}_i^\top {\widehat \Omega}_i^{-1} {\bm{\widehat d}}_i\right) 
\left({\widetilde P}_{\theta_k} {\bm a}_i + \bm{\widetilde \mu}_{\theta_k}\right)^\top {\widetilde\Omega}_i^{-1} {\widetilde P}_i {\bm a}_i\left(-{\bm a}_i^\top {\widetilde P}_i^\top {{\widetilde\Omega}}_i^{-1} {{\widetilde\Omega}}_{\theta_j} {{\widetilde\Omega}}_i^{-1}{\widetilde P}_i {\bm a}_i 
- 2 \bm{\widetilde\mu}_{\theta_k}^\top {\widetilde\Omega}_i^{-1}{\widetilde P}_i {\bm a}_i\right)\Bigg\}.
\end{split}
\end{equation*}
In matrix notation, the derivatives with respect to the sample space used in the Skovgaard's adjustment are 
\begin{equation*}\label{C.2.9}
{\ell'} = 2 R^\top {\bm w}, \:\:\: {U'} = 2\left(R^\top B + V^\top C\right),\:\:\:{\rm and}\:\:\: {\:\:{\J}} = 2\left({\widehat R}^\top F + {\widehat V}^\top G\right),
\end{equation*}
with the elements of $B$, $C$, $F$, $G$, and  ${\bm w}$ being defined in Section  \ref{sec3}.
}
\end{appendix}

\end{document}